\documentclass{article}
\usepackage[utf8]{inputenc}
\usepackage{comment}
\usepackage{changepage}
\usepackage[colorlinks,citecolor=blue,urlcolor=blue,bookmarks=false,hypertexnames=true]{hyperref} 
\title{Long Solutions of Sequence A348480 of the On-Line Encyclopedia of Integer Sequences}
\author{R\"udiger Jehn}
\date{December 2021}

\begin{document}

\maketitle

{\bf Abstract}\\

For numbers $x$ coprime to 10 there exist infinitely many binary numbers $b$ such that the greatest common divisor of $b$ and rev($b$) = $x$ and the sum of digits of $b = x$ (rev($b$) is the digit reversal of $b$). In most cases, the smallest $b$ that fulfill these two constraints contain just a few zeros. But in some cases like for $x = 7, 11, 13$ and 37, $b$ must contain more zeros than ones and these $b$ are called {\it long solutions}. For 11 and 37 it follows directly from the fact that these are porous numbers. For 7 and 13, the proofs that they have long solutions are presented in this paper. 

\section{Introduction}

For numbers $x$ coprime to 10, i.e.\ numbers ending with 1, 3, 7 or 9, there exist infinitely many binary numbers $b$ such that the greatest common divisor of $b$ and rev($b$) = $x$ and the sum of digits of $b = x$. rev($b$) is the digit reversal of $b$, e.g.\ rev(123) = 321. The solutions $b$ converted to decimal numbers form the sequence A348480 of the On-Line Encyclopedia of Integer Sequences \cite{A348480}.\\

For example 
$$gcd(1011,1101)=3$$ 
and since there is no smaller binary number that satisfies this constraint, $b = 1011$ is the solution for $x=3$. $b$ contains 3 ones and 1 zero and is called a short solution because there are more ones than zeros. If there are just a few zeros in the solution a computer search program will find the solution very quickly. However, in the cases $x = 7, 11, 13$ and $37$, there are many more zeros in the solution than ones and only for $x = 7$ and 11, the computer program is able to find the solution in a reasonable amount of time. The solution for $x = 13$ has 73 digits, which is beyond the capabilities of a systematic search.\\

Nevertheless, requirements for the long solutions can be established and the search space will collapse to a single degree of freedom. In the case of the porous numbers $x = 11$ and 37, it was demonstrated in \cite{jehn} that every number $m$ with sum of digits = $x$ and $x$ is a divisor of both $m$ and digit reversal of $m$ must have repeating zeros in their digits. And since a palindrome $m$ where $m$ = rev($m$) like 10101 cannot have 11 or 37 as greatest common divisor of $m$ and rev($m$), more zeros need to be added and the solutions for $x = 11$ and 37 have more zeros than ones.\\

For the sequence A348480, the requirements are even stricter: not only must $x$ be a divisor of both $m$ and digit reversal of $m$, $x$ must also be the greatest common denominator of the two numbers. In the case of 7 and 13, this complicates the search for the smallest solution, because for short numbers $m$ for which $7$ is a divisor of both $m$ and digit reversal of $m$, both $m$ and rev($m$) are also divisible by 13 and vice versa. And therefore no short solution exist for $x = 7$ and 13, as will be shown in this paper.\\

In the last chapter, the special case of $x = 39$ will be treated shortly. 39 does not have a long solution, but it also has no really short solution like all other cases up to and including $x = 71$. The explanation is that the multiplicative order of 10 modulo 39 is 6, the same as the multiplicative order of 10 modulo 7 or modulo 13 and therefore possible solutions are often also divisible by 7 or 13.

\section{Required mathematical tools}
\label{tools}

For our  proof, only a very few simple mathematical tools or fundamentals are required. The first fundamental deals with the parity of numbers. The following statement can be made:\\

If, $x, y, p$ and $q$ are integer numbers and 
\begin{equation}
    x + y = p + q
\label{parity}
\end{equation}
then $x + p - y - q$ is an even number.\\

Proof: \\
$a + b$ has the same parity as $a - b$. Then $x + p - y - q = (x - y) + (p - q)$ is the sum of two numbers with the same parity which gives an even number.\\


A second tool that will be applied are divisibility rules. A number is divisible by 7 if the sum of blocks with length 6 is divisible by 7: You cut the number $b$ in blocks of length 6 starting from the right and sum up the blocks and check if the sum is divisible by 7. The same rule applies for all numbers $x$ where the multiplicative order 10 modulo $x$ is 6, like for 13 and for 39.\\ 

The last required prowess deals with Diophantine equations. For example, if $r$ and $s$ are integers and 
$$5r + 8s \mbox{ is a multiple of } 13$$
then all solutions can be parameterized by:
$$r = 13n + s\mbox{ with } n \mbox{ being any integer.}$$ 

With these tools in hand we can start the proof.

\section{Requirements for a solution for $x = 7$}

The first requirement is that the sum of the digits of the solution $b$ is 7.\\

The second requirement is derived from the fact that 7 divides $b$ since 7 is the greatest common divisor of $x$ and $rev(x)$. The multiplicative order 10 modulo 7 is 6. Hence a number $b$ is divisible by 7 if the sum of blocks with length 6 is divisible by 7.\\ 

Let "$b_{s-1} ... b_3 b_2 b_1 b_0$" be a number $b$ with $s$ digits. We define:\\ 

A = $b_0 + b_6 + b_{12} + $ ... ,

B = $b_1 + b_7 + b_{13} + $ ... , 

C = $b_2 + b_8 + b_{14} + $ ... , 

D = $b_3 + b_9 + b_{15} + $ ... , 

E = $b_4 + b_{10} + b_{16} + $ ... and 

F = $b_5 + b_{11} + b_{17} + $ ...\\

When $10^i$ ($i>1)$ is divided by 7, the first six remainders are 1, 3, 2, 6, 4 and 5 and then this sequence is repeating. We define:
$$
\delta_i = \left\{
    \begin{array}{rl}
        1 & \mbox{if } mod(i,6) = 0  \\
        3 & \mbox{if } mod(i,6) = 1  \\
        2 & \mbox{if } mod(i,6) = 2  \\
        6 & \mbox{if } mod(i,6) = 3  \\
        4 & \mbox{if } mod(i,6) = 4  \\
        5 & \mbox{if } mod(i,6) = 5  \\
    \end{array}
\right.
$$

With these definitions the powers of 10 can be written as:\\

$10^i = j_i \cdot 7 + \delta_i$\\

with some integer $j_i$. Then $b$ can be written as:

$$b = \sum_{i=0}^{s-1} b_i 10^i =  \sum_{i=0}^{s-1} b_i (j_i \cdot 7 + \delta_i) $$\\

Since 7 must divide this sum, the multiples of 7 can be dropped and it follows
\begin{equation}
   \sum_{i=0}^{s-1} b_i \delta_i  = l_1 \cdot 7
   \label{eq:7}
\end{equation}

Hence we get this requirement for the distribution of the 7 ones:

\begin{equation}
    A + 3 B + 2 C + 6 D + 4 E + 5 F = l_1 \cdot 7
   \label{eq:7a}
\end{equation}

The third requirement is derived from the fact that 7 also divides $rev(x)$ since 7 is the greatest common divisor of $x$ and $rev(x)$. As above $rev(b)$ can be written as:

$$b = \sum_{i=0}^{s-1} b_{s-i-1} 10^i =  \sum_{i=0}^{s-1} b_{s-i-1} (j_i \cdot 7 + \delta_i) $$\\

Again, the multiples of 7 can be dropped and it follows
\begin{equation}
   \sum_{i=0}^{s-1} b_{s-i-1} \delta_i  = m_1 \cdot 7
   \label{eq:7rev}
\end{equation}

If 7 divides $b$ then 7 also divides $10b$ and therefore we can add as many zeros to $b$ as needed such that $s$, without loss of generality, can be assumed to be a multiple of 6\footnote{To be on the safe side, the author repeated the proof individually for the other 5 cases of $s$ where $1 \leq mod(s,6) \leq 5$.}. In this case, eq.~\ref{eq:7rev} translates into our third requirement for the distribution of the ones in $b$:

\begin{equation}
    5 A + 4 B + 6 C + 2 D + 3 E + F = m_1 \cdot 7
   \label{eq:7b}
\end{equation}

Adding eq.~\ref{eq:7a} and eq.~\ref{eq:7b} yields:

\begin{equation}
    6 (A + F) + C + D = l_2 \cdot 7
   \label{eq:7plus}
\end{equation}

Note the term $8C$ reduces to $C$ because $7C$ can be subtracted on both sides and is absorbed in a new $l_2$ on the right-hand side of the equation. The same holds for $7E$ and $7B$ which disappear on the left-hand side.\\

Subtracting eq.~\ref{eq:7a} from eq.~\ref{eq:7b} yields:

\begin{equation}
    4 (A + C - D - F) + B - E = m_2 \cdot 7
   \label{eq:7minus}
\end{equation}

In summary, these are the 3 requirements for the positive integers $A, B, C, D, E$ and $F$ such that the greatest common divisor of $b$ and rev($b$) equals 7 and the sum of the digits of $b$ equals 7:

$$
    \begin{array}{rl}
            A + B + C + D + E + F & = 7 \\
                6 (A + F) + C + D & = l_2 \cdot 7 \\
        4 (A + C - D - F) + B - E & = m_2 \cdot 7  \\
    \end{array}
$$

\section{Requirements for a solution for $x = 13$}

The first requirement is that the sum of the digits of the solution $c$ is 13.\\

The second requirement is derived from the fact that 13 divides $c$ since 13 is the greatest common divisor of $x$ and $rev(x)$. The multiplicative order 10 modulo 13 is 6. Hence a number $c$ is divisible by 13 if the sum of blocks with length 6 is divisible by 13.\\ 

Let "$c_{t-1} ... c_3 c_2 c_1 c_0$" be a number $c$ with $t$ digits. We define:\\ 

A' = $c_0 + c_6 + c_{12} + $ ... ,

B' = $c_1 + c_7 + c_{13} + $ ... , 

C' = $c_2 + c_8 + c_{14} + $ ... , 

D' = $c_3 + c_9 + c_{15} + $ ... , 

E' = $c_4 + c_{10} + c_{16} + $ ... and 

F' = $c_5 + c_{11} + c_{17} + $ ...\\

When $10^i$ ($i>1)$ is divided by 13, the first six remainders are 1, 10, 9, 12, 3 and 4 and then this sequence is repeating. We define:
$$
\theta_i = \left\{
    \begin{array}{rl}
        1 & \mbox{if } mod(i,6) = 0  \\
       10 & \mbox{if } mod(i,6) = 1  \\
        9 & \mbox{if } mod(i,6) = 2  \\
       12 & \mbox{if } mod(i,6) = 3  \\
        3 & \mbox{if } mod(i,6) = 4  \\
        4 & \mbox{if } mod(i,6) = 5  \\
    \end{array}
\right.
$$

With these definitions the powers of 10 can be written as:\\

$10^i = k_i \cdot 13 + \theta_i$\\

with some integer $k_i$. Then $c$ can be written as:

$$c = \sum_{i=0}^{t-1} c_i 10^i =  \sum_{i=0}^{t-1} c_i (k_i \cdot 13 + \theta_i) $$\\

Since 13 must divide this sum, the multiples of 13 can be dropped and it follows
\begin{equation}
   \sum_{i=0}^{t-1} c_i \theta_i  = l_3 \cdot 13
   \label{eq:13}
\end{equation}

Hence we get this requirement for the distribution of the 13 ones:

\begin{equation}
    A' + 10 B' + 9 C' + 12 D' + 3 E' + 4 F' = l_3 \cdot 13
   \label{eq:13a}
\end{equation}

The third requirement is derived from the fact that 13 also divides $rev(x)$ since 13 is the greatest common divisor of $x$ and $rev(x)$. As above $rev(c)$ can be written as:

$$c = \sum_{i=0}^{t-1} c_{t-i-1} 10^i =  \sum_{i=0}^{t-1} c_{t-i-1} (k_i \cdot 13 + \theta_i) $$\\

Again, the multiples of 13 can be dropped and it follows
\begin{equation}
   \sum_{i=0}^{t-1} c_{t-i-1} \theta_i  = m_3 \cdot 13
   \label{eq:13rev}
\end{equation}

In this case, again without loss of generality assuming that $t$ is a multiple of 6, eq.~\ref{eq:13rev} translates into our third requirement for $c$:

\begin{equation}
    4 A' + 3 B' + 12 C' + 9 D' + 10 E' + F' = m_3 \cdot 13
   \label{eq:13b}
\end{equation}

Adding eq.~\ref{eq:13a} and eq.~\ref{eq:13b} yields:

\begin{equation}
    5 (A' + F') + 8 (C' + D') = l_4 \cdot 13
   \label{eq:13plus}
\end{equation}

Subtracting eq.~\ref{eq:13a} from eq.~\ref{eq:13b} yields:

\begin{equation}
    3 (A' + C' - D' - F') + 6 (B' - E') = m_4 \cdot 13
   \label{eq:13minus}
\end{equation}

In summary, these are the 3 requirements for the positive integers $A', B', C', D', E'$ and $F'$:

$$
    \begin{array}{rl}
            A' + B' + C' + D' + E' + F' & = 13 \\
            5 (A' + F') + 8 (C' + D')   & = l_4 \cdot 13 \\
    3 (A' + C' - D' - F') + 6 (B' - E') & = m_4 \cdot 13  \\
    \end{array}
$$

\section{Proof that 7 has a long solution}

We will prove that 5 out of the 6 positions $A, B, C, D, E$ and $F$ will be zero and all 7 ones must be located in one position. If we define $R = A + F$ and $S = C + D$ then eq.~\ref{eq:7plus} translates into:
$$ 6 R + S = l_2 \cdot 7$$

This Diophantine equation has the solution space

$$R = 7 n_1 + S$$
with $n_1$ being any integer.\\

Since $0 \leq R, S \leq 7$, $n_1$ can only be -1, 0 or +1. 
\begin{enumerate}
    \item $n_1 = -1 \Rightarrow R = 0$, $S = 7$. This means all ones are in positions $C$ and $D$ and $A = B = E = F = 0$. Then it follows from eq.~\ref{eq:7minus} that $4 (C - D) = m_2 \cdot 7$ or $C - D = m_2 \cdot 7$. Since $C + D = 7$ this is only possible if $C$ or $D$ equals 7. Hence all ones are located in one position, either $C$ or $D$, qed.
    \item $n_1 = +1 \Rightarrow R = 7$, $S = 0$. This means all ones are in positions $A$ and $F$ and $B = C = D = E = 0$. Then it follows from eq.~\ref{eq:7minus} that $4 (A - F) = m_2 \cdot 7$ or $A - F = m_2 \cdot 7$. But this is only possible if $A$ or $F$ equals 7. Hence all ones are located in one position, either $A$ or $F$, qed.
    \item $n_1 = 0 \Rightarrow R = S$. In this case we can conclude two facts: first, if we define $U = A + C - D - F$ then $U$ is even (see Chapter~\ref{tools}) and second, eq.~\ref{eq:13plus} holds because $5 R + 8 S$ is a multiple of 13.
\end{enumerate}

This means either we have proved that 7 has a long solution or $U$ must be even and eq.~\ref{eq:13plus} holds.\\

In the latter case, we define $V = B - E$. Then eq.~\ref{eq:7minus} translates into:
$$ 4 U + V  = m_2 \cdot 7$$

This Diophantine equation has the solution space

$$U = 7 n_2 - 2 V$$
with $n_2$ being any integer.\\

Since $U$ is even, also $n_2$ must be even and since  $-7 \leq U, V \leq 7$, $n_2$ can only be -2, 0 or +2. 
\begin{enumerate}
    \item $n_2 = -2 \Rightarrow U = 0$, $V = -7$. This means all ones are in position $E$, qed.
    \item $n_2 = +2 \Rightarrow U = 0$, $V = 7$. This means all ones are in position $B$, qed.
    \item $n_2 = 0 \Rightarrow U = -2 V$. In this case eq.~\ref{eq:13minus} holds because $3U + 6V = -6V + 6V = 0 = m_4 \cdot 13$.
\end{enumerate}

So either we have a long solution or both equations \ref{eq:13plus} and \ref{eq:13minus} hold. But if both equations are true, then $b$ and rev($b$) are divisible by 13 and the greatest common divisor cannot be 7. Hence 7 has a long solution, qed.\\

The shortest "long solution" for $x=7$ is 

$$1000000000001000001000001000001000001000001,$$ with all ones in position $A$. Digit $b_{36}$ must also be a zero, because the solution cannot be a palindrome, otherwise the greatest common divisor will not be 7.

\section{Proof that 13 has a long solution}

We will prove that 5 out of the 6 positions $A', B', C', D', E'$ and $F'$ will be zero and all 13 ones must be located in one position. If we define $R' = A' + F'$ and $S' = C' + D'$ then eq.~\ref{eq:7plus} translates into:
$$ 5 R' +  8 S' = l_4 \cdot 13$$

This Diophantine equation has the solution space

$$R = 13 n_3 + S$$
with $n_3$ being any integer.\\

Since $0 \leq R, S \leq 13$, $n_3$ can only be -1, 0 or +1. 
\begin{enumerate}
    \item $n_3 = -1 \Rightarrow R' = 0$, $S' = 13$. This means all ones are in positions $C'$ and $D'$ and $A' = B' = E' = F' = 0$. Then it follows from eq.~\ref{eq:13minus} that $3 (C' - D') = m_4 \cdot 13$ or $C' - D' = m_4 \cdot 13$. But this is only possible if $C'$ or $D'$ equals 13. Hence all ones are located in one position, either in $C'$ or $D'$, qed.
    \item $n_3 = +1 \Rightarrow R' = 13$, $S' = 0$. This means all ones are in positions $A'$ and $F'$ and $B' = C' = D' = E' = 0$. Then it follows from eq.~\ref{eq:13minus} that $4 (A' - F') = m_4 \cdot 13$ or $A' - F' = m_4 \cdot 13$. But this is only possible if $A'$ or $F'$ equals 13. Hence all ones are located in one position, either in $A'$ or $F'$, qed.
    \item $n_3 = 0 \Rightarrow R' = S'$. In this case we can conclude two facts: first, if we define $U' = A' + C' - D' - F'$ then $U'$ is even (see Chapter~\ref{tools}) and second, eq.~\ref{eq:7plus} holds because $6 R' + S'$ is a multiple of 7.
\end{enumerate}

This means either we have proved that 13 has a long solution or $U'$ must be even and eq.~\ref{eq:7plus} holds.\\

In the latter case, we define $V' = B' - E'$. Then eq.~\ref{eq:13minus} translates into:
$$ 3 U' + 6 V'  = m_4 \cdot 13$$

This Diophantine equation has the solution space

$$U' = 13 n_4 - 2 V'$$
with $n_4$ being any integer.\\

Since $U'$ is even, also $n_4$ must be even and since  $-13 \leq U', V' \leq 13$, $n_4$ can only be -2, 0 or +2. 
\begin{enumerate}
    \item $n_4 = -2 \Rightarrow U' = 0$, $V' = -13$. This means all ones are in position $E'$, qed.
    \item $n_4 = +2 \Rightarrow U' = 0$, $V' = 13$. This means all ones are in position $B'$, qed.
    \item $n_4 = 0 \Rightarrow U' = -2 V'$. In this case eq.~\ref{eq:7minus} holds because $4U' + V' = -8V' + V' = -7V'$ is a multiple of 7.
\end{enumerate}

So either we have a long solution or both equations \ref{eq:7plus} and \ref{eq:7minus} hold. But if both equations are true, then $b$ and rev($b$) are divisible by 7 and the greatest common divisor cannot be 13. Hence 13 has a long solution, qed.

\section{39 has quite a long "short solution"}

Since the multiplicative order 10 modulo 39 is also 6, it turns out that for most numbers $b$ where the sum of digits of $b$ is 39 and 39 divides both $b$ and rev($b$), also 7 divides both $b$ and rev($b$) and therefore 39 is not the greatest common divisor. In fact the smallest $b$ where the greatest common divisor is 39 has 16 zeros, i.e.\ $b$ has 55 digits in total. With a computer program, an iterative search was performed given these three constraints:

$$
    \begin{array}{rl}
        A + B + C + D + E + F & = 39 \\
        A + 10 B + 22 C + 25 D + 16 E + 4 F & = l_5 \cdot 39 \\
        4 A + 16 B + 25 C + 22 D + 10 E + F & = m_5 \cdot 39  \\
    \end{array}
$$
The result was $A = 10$, $B = C = 9$, $D = 6$, $E = 0$ and $F = 5$. In binary terms, this translates into 
b = 1000111000111000111001111101111101111101111101111101111 or in decimal form 20016007615544303 which is term 16 of A348480.

\bibliography{references} 
\bibliographystyle{ieeetr}

\end{document}